\documentclass[11pt]{amsart}

\usepackage{amssymb,amsmath}
\usepackage{bbm}
\usepackage{graphics}
\usepackage{a4wide}

\newtheorem{theorem}{Theorem}
\newtheorem{prop}{Proposition}
\newtheorem{lemma}{Lemma}
\newtheorem{coro}{Corollary}
\newtheorem{fact}{Fact}
 
{
\theoremstyle{definition}

\newtheorem{remark}{Remark}

}

\newcommand{\ts}{\hspace{0.5pt}}

\newcommand{\RR}{\mathbb{R}\ts}
\newcommand{\QQ}{\mathbb{Q}\ts\ts}
\newcommand{\ZZ}{\mathbb{Z}}
\newcommand{\NN}{\mathbb{N}}
\newcommand{\HH}{\mathbb{H}\ts}
\newcommand{\II}{\mathbb{I}}

\newcommand{\LL}{\mathcal{L}}
\newcommand{\alat}{L}
\newcommand{\up}{\widetilde{\!\!\!\hphantom{m}}}
\newcommand{\uu}{{\ts\ts\!\times\!}}
\newcommand{\bt}{{\scriptscriptstyle \bullet}}

\newcommand{\ii}{\mathrm{i}}
\newcommand{\jj}{\mathrm{j}}
\newcommand{\kk}{\mathrm{k}}

\newcommand{\gG}{\varGamma}
\newcommand{\gL}{\varLambda}

\newcommand{\op}{{\dot{+}}}

\DeclareMathOperator{\lcm}{lcm}
\DeclareMathOperator{\tr}{tr}
\DeclareMathOperator{\nr}{nr}
\DeclareMathOperator{\N}{N}
\DeclareMathOperator{\Tr}{Tr}

\begin{document}

\title[Similar sublattices of $A_4$]
{Similar sublattices of the root lattice $A_4$}

\author{Michael Baake}
\author{Manuela Heuer}
\author{Robert V. Moody}

\address{Fakult\"at f\"ur Mathematik, Universit\"at Bielefeld, \newline
\hspace*{12pt}Box 100131, 33501 Bielefeld, Germany}
\email{mbaake@math.uni-bielefeld.de}

\address{Department of Mathematics and Statistics, The Open University, 
\newline
\hspace*{12pt}Walton Hall, Milton Keynes MK7 6AA, UK}
\email{m.heuer@open.ac.uk}

\address{Department of Mathematics and Statistics, 
University of Victoria, \newline
\hspace*{12pt}Victoria, British Columbia V8W 3P4, Canada}
\email{rmoody@uvic.ca}

\begin{abstract} 
  Similar sublattices of the root lattice $A_4$ are possible
  \cite{CRS} for each index that is the square of a non-zero integer
  of the form $m^2 + mn - n^2$. Here, we add a constructive approach,
  based on the arithmetic of the quaternion algebra $\HH (\QQ
  (\sqrt{5}))$ and the existence of a particular involution of the
  second kind, which also provides the actual sublattices and the
  number of different solutions for a given index. The corresponding
  Dirichlet series generating function is closely related to the zeta
  function of the icosian ring.
\end{abstract}

\maketitle

\section{Introduction}

Any lattice $\gG$ in Euclidean space $\RR^d$ possesses sublattices
that are images of $\gG$ under a non-zero similarity, i.e., under a
linear map $S$ of $\RR^d$ with $\langle S u | S v \rangle = c \,
\langle u | v\rangle$ for all $u,v\in\RR^d$, where $c>0$ and $\langle
. | .  \rangle$ denotes the standard Euclidean scalar product.  If
$\gG'=S\gG\subset\gG$, $\gG'$ is a \emph{similar sublattice} (or
similarity sublattice, SSL) of $\gG$, with index $[\gG : \gG'] =
c^{\ts d/2}$. If the lattice $\gG$ has a rich point symmetry
structure, there are usually interesting SSLs, beyond the trivial ones
of the form $m\gG$ with $m\in\NN$.  Several examples have been
investigated in the past, compare \cite{BM1,BG,CRS,BM2} and references
given there. Of particular interest are root lattices, due to their
interesting geometry and their abundance in mathematics and its
applications. In \cite{CRS}, the possible values of $c$ were
classified for many root lattices, while more complete information in
terms of generating functions is available in dimension $d\le 4$, see
\cite{BM1,BM2}.

One important example that has not yet been solved completely is the
root lattice $A_4$. It is the purpose of this short article to close
this gap. After a short summary of a more general nature, we review
the mathematical setting in $4$-space via quaternions.  This permits
the solution of the problem in Section~\ref{solution}, using the
arithmetic of a maximal order of the division algebra $\HH
(\QQ(\sqrt{5}\,))$.  Related results will briefly be summarized
afterwards, followed by some details on the underlying algebraic
structure that we found along the way.

\section{Generalities}

A lattice $\gG$ in Euclidean space $\RR^d$ is a free $\ZZ$-module of
rank $d$ whose $\RR$-span is all of $\RR^d$. Two lattices in $\RR^d$
are called \emph{similar} if they are related by a non-zero
similarity. For a given lattice $\gG\subset\RR^d$, the similar
sublattices (SSLs) form an important class, and one is interested in
the possible indices and the arithmetic function $g^{}_{\gG}(m)$ that
counts the distinct SSLs of index $m$. Clearly,  $g^{}_{\gG}(1)=1$. 
One can show \cite{H} that
$g^{}_{\gG}(m)$ is super-multiplicative, i.e., $g^{}_{\gG}(mn)\ge
g^{}_{\gG}(m) g^{}_{\gG}(n)$ for coprime $m,n$.  It is not always
multiplicative, as one can see from considering suitable rectangular
lattices in the plane \cite{Z}, such as the one spanned by $(2,0)$ and
$(0,3)$. However, in many relevant cases, and in all cases that will
appear below, $g^{}_{\gG}(m)$ is multiplicative.  This motivates the
use of a Dirichlet series generating function for it,
\begin{equation}
   D_{\gG} (s) \; := \; \sum_{m=1}^{\infty}
   \frac{g^{}_{\gG}(m)}{m^s} \, ,
\end{equation} 
which will then have an Euler product expansion as well.

If $\gG$ is a generic lattice in $\RR^d$, the only SSLs of $\gG$ will
be its integrally scaled versions, i.e., the sublattices of the form
$m\gG$ with $m\in\NN$. In this case, as $m\gG$ has index $m^d$ in
$\gG$, the generating function simply reads
\begin{equation} \label{generic}
   D_{\gG} (s) \; = \; \zeta(ds) \, ,
\end{equation}
where $\zeta(s)$ is Riemann's zeta function. If SSLs other than
those of the form $m\gG$ exist, each of them can again be scaled
by an arbitrary natural number, so that 
\begin{equation} 
   D^{}_{\gG} (s) \; = \; \zeta(ds) \, D^{\sf pr}_{\gG} (s) \, ,
\end{equation}
where $D^{\sf pr}_{\gG} (s)$ is the Dirichlet series generating 
function for the \emph{primitive} SSLs of $\gG$ (we shall define this 
in more detail below).

Let us first summarize some general properties of the generating
functions in this context.  The following statement follows from a
simple conjugation argument.

\begin{fact} \label{invar-similar}
If $\gG$ and $\gL$ are similar lattices in\/ $\RR^d$,
one has $D_{\gG} (s) = D_{\gL} (s)$.  \qed
\end{fact}

Given a lattice $\gG\subset\RR^d$, its \emph{dual} lattice $\gG^*$ is
defined as
\begin{equation} \label{dual}
    \gG^* \; := \; \{ x\in\RR^d \mid \langle x|y\rangle\in\ZZ
    \mbox{ for all } y\in\gG\}\, .
\end{equation}

\begin{fact} \label{invar-dual}
If $\gG$ is a lattice in\/ $\RR^d$, one has $D_{\gG} (s) = D_{\gG^*} (s)$.
\end{fact}
\begin{proof}
If $S$ is any (non-zero) linear similarity in $\RR^d$, one has
\[
   S\gG\subset\gG \quad \Longleftrightarrow \quad
   S^t\gG^*\subset\gG^* ,
\]
which, in view of \eqref{dual}, is immediate from the relation 
$\langle x|Sy\rangle = \langle S^t x|y\rangle$.

Observing $\det(S)=\det(S^t)$, one thus obtains an index preserving
bijection between the SSLs of $\gG$ and those of $\gG^*$, whence we
have $g^{}_{\gG} (m) = g^{}_{\gG^*} (m)$ for all $m\in\NN$. This gives
$D_{\gG} (s) = D_{\gG^*} (s)$.
\end{proof}

\begin{remark} \label{duality}
  Defining the transposed similarity via the scalar product makes it
  transparent that there is a certain scaling degree of freedom to
  choose the non-degenerate bilinear form without changing the result
  of Fact~\ref{invar-dual}, which is in line with
  Fact~\ref{invar-similar}.
\end{remark}

\section{Mathematical setting in $4$-space}

\begin{figure}
 \scalebox{0.6}{\includegraphics{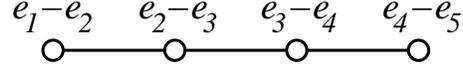}}
 \caption{Standard basis representation of the root lattice $A_4$.
 \label{a4-fig}}
 \end{figure}

As is apparent from the treatment in \cite{CRS}, the usual description
of $A_4$ as a lattice in a $4$-dimensional hyperplane of $\RR^5$, see
Figure~\ref{a4-fig}, gives access to the possible values of $c$, and
hence to the possible indices, but not necessarily to the actual
similarities and the SSLs.  Therefore, we use a description in $\RR^4$
instead, which will enable us to use quaternions, based on the inclusion
of a root system of type $A_4$ within one of type $H_4$, see
\cite[Prop.~2]{S1} or \cite[Thm.~4.1]{CMP}. In particular, we
shall need the Hamiltonian quaternion algebras $\HH (\RR)$, $\HH
(\QQ)$ and $\HH (K)$, with the quadratic field $K = \QQ(\sqrt{5}\,)$,
generated by the basis quaternions $1,\ii,\jj,\kk$.  The latter
observe Hamilton's relations
\[
    \ii^2 \; = \ \jj^2 \; = \; \kk^2 \; = \; \ii\jj\kk \; = \; -1
\]
and will be identified with the canonical Euclidean basis in
$4$-space.  Quaternions are thus also written as vectors,
$x=(x^{}_{0},x^{}_{1},x^{}_{2},x^{}_{3})$, and the algebra is equipped
with the usual conjugation $x\mapsto\bar{x}$, where
$\bar{x}=(x^{}_{0},-x^{}_{1},-x^{}_{2},-x^{}_{3})$. For later use, we
also introduce the reduced trace $\tr (x)$ and the reduced norm $\nr
(x)$ in a quaternion algebra via
\[
     \tr (x) \, := \, x + \bar{x} \, , \quad
     \nr (x) \, := \, x \bar{x} = \lvert x \rvert^2 .
\]
If we work over the field $K=\QQ(\sqrt{5}\,)$, we shall also need
its trace and norm, defined as
\[
     \Tr (\alpha) \, := \, \alpha + \alpha'  \, , \quad
     \N (\alpha) \, := \, \alpha\ts\alpha' ,
\]
where ${}^{\prime}$ is the algebraic conjugation in $K$, defined
by $\sqrt{5}\mapsto - \sqrt{5}$.

The first benefit of this setting is the following result.
\begin{fact} \label{quat-simi} 
  All non-zero linear similarities in\/ $\RR^4 \simeq \HH (\RR)$
  can be written as either $x\mapsto p x q$ $($orientation
  preserving case\/$)$ or as $x\mapsto p \bar{x} q$ $($orientation
  reversing case\/$)$, with non-zero quaternions $p,q \in \HH (\RR)$.
  The determinants of the linear maps defined this way are
  $\pm \lvert p\rvert^4 \lvert q\rvert^4$. Conversely, all mappings
  of the form $x\mapsto p x q$ and $x\mapsto p \bar{x} q$ are
  linear similarities.
\end{fact}

\begin{proof}
This follows easily from the corresponding result on orthogonal
matrices in $\RR^4$, as every non-zero similarity is the product
of an orthogonal transformation and a non-zero homothety, see
\cite{KR,BM2} for details. The result on the determinants is standard.
\end{proof}

For some purposes, it is convenient to refer to the standard matrix
representation of the linear map $x\mapsto pxq$, as defined via  
$(p x q)^t  =  M(p,q)\ts x^t$. Details can be found in \cite{KR,BM2}.

Following \cite{CMP}, we consider the lattice
\begin{equation} \label{lat-1}
   \alat \; = \; \bigl\langle (1,0,0,0),
   \tfrac{1}{2}(-1,1,1,1),(0,-1,0,0),\tfrac{1}{2}
   (0,1,\tau \! - \! 1,-\tau)\bigr\rangle_{\ZZ}
\end{equation}
where $\tau=(1+\sqrt{5}\,)/2$ is the golden ratio. Due to $\Tr
(\tau)=1$ and $\N (\tau) = -1$, $\tau$ is an integer in $K$. In fact,
$\tau$ is a fundamental unit of $\ZZ[\tau]$, the ring of
integers of $K$.

The Gram matrix of $\alat$, calculated with $\langle.|.\rangle$,  reads
\begin{equation} \label{gram-matrix}
  \frac{1}{2}     \begin{pmatrix}
     2 & -1 & 0 & 0 \\ -1 & 2 & -1 & 0 \\
     0 & -1 & 2 & -1 \\ 0 & 0 & -1 & 2
     \end{pmatrix} ,
\end{equation}
which shows that $\alat$ is a scaled copy of the root lattice $A_4$ in
its standard representation with basis vectors of squared length $2$,
compare \cite{CS} and Figure~\ref{a4-fig}.  Our choice of this
particular realization of the $A_4$ root lattice is motivated by the
observation that $\alat$ is a subset of the so-called \emph{icosian
  ring} $\ts\II$, which has a powerful arithmetic structure. A
concrete representation of $\II$ is obtained via the vectors
\begin{equation} \label{h4-vectors}
   (\pm 1,0,0,0) \, , \;
   \tfrac{1}{2} (\pm 1,\pm 1,\pm 1,\pm 1) \, , \;
   \tfrac{1}{2} (\pm \tau',\pm \tau,0,\pm 1) 
\end{equation}
and all their \emph{even} coordinate permutations. Together, they
form a set of $120$ vectors (all of length $1$) which constitute a
non-crystallographic \emph{root system of type $H_4$},
denoted by $\varDelta_{H_4}$, see \cite{CMP,BM2} for details. 
The orthogonal transformations $\{ x\mapsto pxq\, , \; x\mapsto
p\bar{x}q \mid p,q \in \varDelta_{H_4}\}$ form the corresponding
\emph{Coxeter group $W_{H_4}$} of all $(120)^2$ point
symmetries of $\varDelta_{H_4}$. 

The icosian ring $\II$ is the $\ZZ$-span of $\varDelta_{H_4}$, 
$\II := \langle \varDelta_{H_4}\rangle_{\ZZ}$, and has rank $8$.
At the same time, $\II$  is  a $\ZZ[\tau]$-module of rank 
$4$, and can alternatively be written as
\begin{equation} \label{icosian-1}
   \II \; = \; \bigl\langle (1,0,0,0),(0,1,0,0),
   \tfrac{1}{2}(1,1,1,1),\tfrac{1}{2}(1 \! - \!\tau,\tau,0,1)
   \bigr\rangle_{\ZZ[\tau]} \, .
\end{equation}
In this formulation, $\II$ is a maximal order in the quaternion algebra 
$\HH (K)$, which has class number $1$, see \cite{V}. Consequently, 
all ideals of $\II$ are principal, and one has a powerful notion of
prime factorisation.

\begin{remark} \label{rem-one} 
  The quadratic form $\tr (x\bar{y}) = 2\ts\langle x|y\rangle$ takes
  only integer values on the icosian ring, and is thus the canonical
  form on $\II$ in the context of root systems. Relative to this form,
  $\alat$ \emph{is} the root lattice $A_4$, with the weight lattice
  $A_4^*$ as its dual (defined relative to the trace form) and the index
  relation $[A_4^* : A_4^{} ] = 5$, see \cite{CS,CMP} for details.
\end{remark}

The icosian ring $\II$ also contains the $\ZZ[\tau]$-module
\begin{equation} \label{def-L}
    \LL \; := \; \bigl\langle 1,\ii,\jj,\kk \bigr\rangle_{\ZZ[\tau]}
    \; \simeq \; \ZZ[\tau]^4
\end{equation}
as a submodule of index $16$. The same statement applies to $\II'$,
which is another maximal order in $\HH(K)$. Note that $\LL$ is
an index $4$ submodule of $\II\cap\II'$, the latter thus being of
index $4$ both in $\II$ and in $\II'$.

The benefit of this approach will be that we can use the arithmetic of
the ring $\II$.  Let us recall the following result from \cite{BM2},
which is also immediate from Eqs.~\eqref{icosian-1} and \eqref{def-L}.
\begin{fact} \label{L-in-I}
   $2\ts\II\subset\LL$, i.e., $x\in\II$ implies $2x\in\LL$. \qed
\end{fact}

The golden key for solving the SSL problem is the mapping
$\,\widetilde{.}  : \, \HH (K) \longrightarrow \HH (K)$,
$x\mapsto\widetilde{x}$ (also denoted by $\eta$ later on), with
\begin{equation} \label{the-map}
   \widetilde{x} \; := \; (x^{\ts\prime}_{0},x^{\ts\prime}_{1},
     x^{\ts\prime}_{3},x^{\ts\prime}_{2}) \, .
\end{equation}
The relevance of this map was noticed in \cite{H} as a result of
explicit calculations around SSLs with small indices.  The rather
strange looking combination of a permutation of two coordinates with
algebraic conjugation of all coordinates, which we shall call a
\emph{twist map} from now on, is an algebra involution of the second
kind \cite[Ch.~I.2]{Rost}. It has the following important properties
(see also Section~\ref{after} for more).
\begin{lemma} \label{props}
  The twist map\/ $\widetilde{.}$ of $\eqref{the-map}$ is a\/
  $\QQ$-linear and $K$-semilinear involutory algebra anti-automorphism,
  i.e., for arbitrary $x,y\in\HH (K)$ and
  $\alpha\in K$, it satisfies:\
\begin{itemize}
\item[\rm (a)]  $\widetilde{x+y}=\widetilde{x}+\widetilde{y}\;$ 
   and $\;\widetilde{\alpha x} = \alpha'\ts\widetilde{x}$,
\item[\rm (b)] $\widetilde{xy} = \widetilde{y}\,\widetilde{x} \;$
   and $\; \widetilde{\widetilde{x}} = x$,
\item[\rm (c)] $\widetilde{\bar{x}} = \overline{\widetilde{x}}$ and
   thus, for $x\neq 0$, also $(\widetilde{x}\ts)^{-1} =  
   \widetilde{x^{-1}}$.
\end{itemize}
It maps $K$, the centre of the algebra\/ $\HH (K)$, onto itself, but
fixes only the elements of\/ $\QQ$ within $K$.  Moreover, the icosian
ring is mapped onto itself, $\widetilde{\II}=\II$, as is its algebraic
conjugate, $\II'$.
\end{lemma}
\begin{proof}
Most of these properties are immediate from the definition, and
(b) can be proved by checking the action of the twist map on
the basis quaternions $1,\ii,\jj,\kk$. The statements on $K$, $\II$
and $\II'$ follow easily from the definition in \eqref{icosian-1}.
\end{proof}

Seen as a vector space over $\QQ$, $\HH(K)$ has dimension $8$ and
can be split into a direct sum, $\HH(K) = V_+ \op V_-$, where
\[
    V_{\pm} \; := \; \{ x\in\HH (K) \mid \widetilde{x} = \pm x \}
\]
are the eigenspaces under the twist map, with $V_{+}\cap V_{-} =
\{ 0 \}$ and $V_{-} = \sqrt{5}\, V_{+}$. Concretely, one can choose
the $4$ basis vectors of $\alat$ from \eqref{lat-1} also as a
$\QQ$-basis for $V_{+}$. 

\begin{remark}
  Defining $u=\sqrt{5}\, \ii$, $v=\jj - \kk$ and $w=\sqrt{5}\ts (\jj +
  \kk)$, one checks that $u^2 = -5$, $v^2 = -2$ and $w = uv = - vu$,
  so that $D:=\langle 1,u,v,w\rangle_{\QQ}$ is a quaternion algebra
  over $\QQ$, of type $H (\frac{-5,-2}{\QQ})$ in the terminology of
  \cite[Sec.~57]{O}. Extending this algebra by adjoining $\sqrt{5}$
  leads to $\HH(K)$, so $\HH(K) = \QQ(\sqrt{5}\,)\otimes^{}_{\QQ} D$.
  The benefit \cite{Rehmann} is that the involution
  $\,\widetilde{.}\,$ looks more canonical in this setting. In fact,
  $\widetilde{u} = -u = \bar{u}$, and analogously for $v$ and $w$.
  Consequently, for $\alpha \in K$ and $x\in D$, one has $(\alpha
  x)^{\up} = \alpha'\ts\bar{x}$, which clearly shows the nature of the
  involution. Note, however, that $\alat$ is a subset of $V_{+}$ and
  is thus not contained in $D$.
\end{remark}

These observations, together with the nature of $\II$ as a maximal
order in $\HH(K)$, suggest that $\alat$ might actually be the 
$\ZZ$-module of fixed points of the twist map inside $\II$. 
\begin{prop} \label{key-prop}
   One has\/ $\II \cap V_+  \, = \,
   \{ x\in\II\mid \widetilde{x}=x \}\, = \, \alat $, with $\alat$
   as in $\eqref{lat-1}$.
\end{prop}
\begin{proof}
The basis vectors of $\alat$ in \eqref{lat-1} are fixed under the
twist map, so $\widetilde{x}=x$ is clear for all $x\in\alat$. We also
know from Lemma~\ref{props} that $\widetilde{\II}=\II$. To prove
our claim, we have to show that no element of $\II\setminus\alat$
is fixed under the twist map.

Let $\alat [\tau]$ be the $\ZZ[\tau]$-span of $\alat$. It is clear that
\[
   \alat [\tau] \; = \; \alat + \tau\alat 
   \; = \; \alat \,\op \,\tau\alat \, ,
\]
where the latter equality follows from $\alat\cap\tau\alat=\{0\}$, and
also that $\widetilde{\alat[\tau]}=\alat[\tau]$.  Observe that
$\lvert\det_{K} (\II)\rvert = 1/4$ and $\lvert\det_{K}
(\alat[\tau])\rvert = \sqrt{5}/4$, where the $K$-determinants are just
the ordinary determinants of the bases given in \eqref{lat-1} and
\eqref{icosian-1}. Using the connection of the corresponding index
with the ratio of the $\QQ$-determinants, see \cite{BPR} for details,
one concludes that
\begin{equation} \label{taulat-in-ico}
     \bigl[\II : \alat [\tau] \bigr] \; =\; 
     \bigl| \N \bigl( \sqrt{5}\,\bigr) \bigr|
     \; = \; 5 \, ,
\end{equation}
so that $\II/\alat [\tau]\simeq C_5$. 

Consequently, the twist map, which is an involution, induces an
automorphism on the cyclic group $C_5$. The order of this automorphism
must divide $2$.  Since $0\neq v=\jj - \kk \in \II
\setminus\alat[\tau]$, but satisfies $\widetilde{v}=-v$, the induced
automorphism on $C_5$ cannot be the identity. This leaves only
inversion (i.e., $k\mapsto -k$ mod $5$), which has no fixed point in
$C_5$ other than $0$, so that all fixed points of the twist map inside
$\II$ must lie in $\alat[\tau]$. It is easy to see that a quaternion
of the form $x+\tau y$ with $x,y\in\alat$ is fixed if and only if
$y=0$, which completes the argument.
\end{proof}

\begin{remark} \label{remark-dual} An alternative way to prove
  Proposition~\ref{key-prop} would be to use Remark~\ref{rem-one} and
  to show that all potential $x\in\II$ with $\widetilde{x}=x$ are
  elements of the dual of $\alat$ (taken with respect to $\tr
  (x\bar{y})$), but that the coset representative of $A_4^*/A_4^{}$
  cannot be chosen in $\II$. In fact, extending $A_{4}^{}$ and
  $A_{4}^{*}$ to lattices over $\ZZ[\tau]$, compare \cite{S2}, one
  sees that $\II$ is intermediate between them, but is still integral
  for $\tr (x\bar{y})$. When viewed as a $\ZZ[\tau]$-module, $\II$ has
  class number $1$, see \cite[Thm.~3.4]{S2}. Moreover, it is a
  principal ideal domain and a maximal order in $\HH(K)$. The intimate
  relation with quaternion arithmetic is the entry point for our
  further analysis.
\end{remark}

In fact, in the setting of Remarks \ref{rem-one} and
\ref{remark-dual}, we have $\alat = A_4$ and thus $\alat[\tau] =
\ZZ[\tau]\otimes_{\ZZ} A_4$. Then, its dual over $\ZZ[\tau]$ (relative
to $\tr (x\bar{y})$ again) is $\ZZ[\tau]\otimes_{\ZZ} A_4^*$, which
has $K$-index $5$ over $\alat[\tau]$, where the $K$-index (up to units
in $\ZZ[\tau]$) is given by the determinant of the linear mapping that
sends a basis of $\alat[\tau]$ to one of $\alat[\tau]^*$. The
various inclusions can now be summarized as
\[
   2\alat[\tau] \; \stackrel{16}{\subset} \; \alat[\tau] \;
   \stackrel{\sqrt{5}}{\subset} \; \II \; 
   \stackrel{\sqrt{5}}{\subset} \; \alat[\tau]^*
\]
together with
\[   
    2\alat[\tau] \; \stackrel{\sqrt{5}}{\subset} \; 2 \II \;
    \stackrel{4}{\subset} \; \LL \; \stackrel{4}{\subset} 
    \; \II \;  \stackrel{\sqrt{5}}{\subset} \; \alat[\tau]^* ,
\]
where the integers on top of the inclusion symbols denote the
corresponding $K$-indices. The indices over $\QQ$ are then
just the squares of these numbers.

Also, as a result of Proposition~\ref{key-prop}, we get that
\begin{equation} \label{root-systems}
      \varDelta_{A_4} \; = \; V_{+} \cap \varDelta_{H_4} 
\end{equation}
with the root system $\varDelta_{H_4}$ as described above.  Explicitly,
observing $\tau'=1-\tau$, this leaves us with the $20$ roots of the root
system $\varDelta_{A_4}$,
\begin{equation} \label{a4-roots}
  \pm (1,0,0,0) , \, \pm (0,1,0,0) , \,
  \pm \tfrac{1}{2}\ts (\pm 1,\pm 1,1,1)  , \,
  \pm \tfrac{1}{2}\ts (0,\pm 1,\tau\!-\! 1,-\tau) , \,
  \pm \tfrac{1}{2}\ts (\pm 1,0,-\tau,\tau\!-\! 1) ,
\end{equation}
which will become important later on, see \cite[Prop.~2]{S1}
and \cite{CMP} for the relation between the corresponding
reflection groups.

\section{Similarities for $A_4$ via quaternions}
\label{solution}

We are interested in the SSLs of the root lattice $A_4$. They are in
one-to-one relation to those of the lattice $\alat$ defined in
Eq.~\eqref{lat-1}. For convenience, we define
$\HH (K)^\bt_{} = \HH (K)\setminus \{0\}$, and analogously for
other rings.

\begin{lemma} \label{para-1}
   All SSLs of the lattice $\alat$, defined as in Eq.~$\eqref{lat-1}$,
   are images of $\alat$ under orientation preserving mappings
   of the form $x\mapsto pxq$, with $p,q \in
   \HH (K)^\bt$ and $K=\QQ(\sqrt{5}\,)$.
\end{lemma}
\begin{proof}
As $\alat$ is invariant under the conjugation map $x\mapsto\bar{x}$,
we need not consider orientation reversing similarities. By
Fact~\ref{quat-simi}, all SSLs are thus images of $\alat$ under
mappings $x\mapsto pxq$ with $p,q\in\HH (\RR)^\bt$.

A simple argument with the matrix entries of the linear
transformations defined by the mappings $x\mapsto px q$, compare
\cite{BM2}, shows that all products $p_i q_j$ must be in the quadratic
field $K=\QQ(\sqrt{5}\,)$, which leaves the choice to take $p,q\in \HH
(K)^\bt$.
\end{proof}

We need to select an appropriate subset of mappings that reach
all SSLs of $\alat$. A first step is provided by the following
observation.

\begin{lemma} \label{all-do}
    If $p\in\II$, $p\alat \widetilde{p}$ is an SSL of $\alat$.
\end{lemma}

\begin{proof}
By Fact~\ref{quat-simi}, $p\alat\widetilde{p}$ is similar to $\alat$,
so it remains to be shown that $p\alat\widetilde{p}\subset\alat$.
Observe that $p\in\II$ implies $\widetilde{p}\in\II$, so that $p\alat
\widetilde{p} \subset\II$ is clear. If $x$ is any point of $\alat$, we
have $\widetilde{x}=x$ by Proposition~\ref{key-prop}.  Using the
properties of the twist map from Lemma~\ref{props}, one gets
\[
   (p x \widetilde{p}\, )^{\up} \, = \,
   \widetilde{\widetilde{p}}\, \widetilde{x}\, \widetilde{p} 
   \, = \, p\ts x \widetilde{p} \, .
\]
Consequently, again by Proposition~\ref{key-prop}, $p x \widetilde{p}
\in \alat$, and hence $p\alat \widetilde{p} \subset \alat$, as
claimed.
\end{proof}

\begin{prop}  \label{fundamental-para}
   If $p \alat q \subset \alat$ with $p,q\in\HH (K)^\bt$,
   there is an $\alpha\in\QQ$ such that
\[
   q \; = \; \alpha\, \widetilde{p} \, .
\]
\end{prop}
\begin{proof}
If $p,q\in\HH(K)^\bt$, the inclusion $p \alat q\subset\alat$ implies
$pxq = \widetilde{pxq} = \widetilde{q} x \widetilde{p}$ for all 
$x\in\alat$, hence also $ x \widetilde{p} q^{-1}
= (\widetilde{q}\ts )^{-1}p x = (\widetilde{p} q^{-1})^{\up} x$.
Since $1\in\alat$, this implies
$\widetilde{p} q^{-1} = (\widetilde{p} q^{-1})^{\up}$,
and we get
\[
   x \widetilde{p} q^{-1} \; = \; \widetilde{p} q^{-1} x \, ,
\]
still for all $x\in\alat$. Noting that $\langle\alat\rangle_K=\HH (K)$,
the previous equation implies that $\widetilde{p} q^{-1}$ must be
central, i.e., an element of $K$. Consequently, $q=\alpha
\widetilde{p}$ for some $\alpha\in K$.

Since $p\in\HH (K)$, we can choose some $0\neq\beta\in\ZZ[\tau]$ such
that $w=\beta p\in\II$. Observing that $(\beta p)^{\up}=
\beta'\ts\widetilde{p}$, one sees that $\alpha \, p x\widetilde{p} =
\frac{\alpha}{\N (\beta)}\, w x \widetilde{w}$, where $0\neq \N
(\beta)\in\ZZ$. As $w\alat\widetilde{w} \subset\alat$ by
Lemma~\ref{all-do}, and since $\alat \cap \tau \alat = \{0\}$, 
the original relation $p\alat q\subset \alat$ now implies
$\frac{\alpha}{\N (\beta)}\in\QQ$, hence also $\alpha\in\QQ$.
\end{proof}

An element $q\in\II$ is called \emph{$\II$-primitive} when all
$\ZZ[\tau]$-divisors of $q$ are units. Another way to phrase this is
to say that the $\II$-content of $q$,
\begin{equation} \label{I-cont}
    \mathrm{cont}^{}_{\II} (q) \; := \; \mathrm{lcm}
    \{ \alpha\in\ZZ[\tau]^\bt \mid q \in \alpha \II \}\, ,
\end{equation}
is a unit in $\ZZ[\tau]$. Whenever the context is clear, we shall
simply speak of primitivity and content of $q$. Note that the content
is unique up to units of $\ZZ[\tau]$, as $\ZZ[\tau]$ is Euclidean and
hence also a principal ideal domain. More generally, one has to formulate
the content via a fractional ideal, see \cite{BM2} for details.

\begin{coro} \label{para-2}
  All SSLs of the lattice $\alat$ are images of $\alat$ under mappings
  of the form $x\mapsto \alpha \ts p x \widetilde{p}$ with $p \in\II$
  primitive and $\alpha\in\QQ$.
\end{coro}

\begin{proof}
  By Proposition~\ref{fundamental-para}, we know that mappings of the
  form $x\mapsto \beta q x \widetilde{q}$ with $q\in\HH (K)$ and
  $\beta\in\QQ$ suffice.  Since all coordinates of $q$ are in
  $K=\QQ(\sqrt{5}\,)$, there is a natural number $m$ such that $p:=mq
  \in \LL$, with $\LL$ as in \eqref{def-L}. Then, with $\alpha: =
  \beta/m^2$, one has $\alpha \ts px \widetilde{p} = \beta\ts qx
  \widetilde{q}$ whence the mappings $x\mapsto\alpha \ts px
  \widetilde{p}$ and $x\mapsto \beta\ts qx \widetilde{q}$ are equal.
  While $\alpha$ is still in $\QQ$, $p$ and $\widetilde{p}$ are now
  elements of $\LL\subset\II$.

  If $p$ is primitive in $\II$, we are done. If not, we know that
  $p/c$ is a primitive element of $\II$ when $c=\mathrm{cont}^{}_{\II}
  (p)$.  Simultaneously, we have $(p/c)^{\widetilde{\hphantom{a}}} =
  \widetilde{p}/{c\ts}'$. Since $c\in\ZZ[\tau]$, we know that $cc' \in
  \ZZ$, so that this factor can be absorbed into $\alpha$.
\end{proof}

At this stage, we recollect an important property of the icosian
ring from \cite{BM2}.
\begin{fact} \label{no-escape}
  Let $p,q\in\HH (K)^\bt$ be such that $p\ts\II q\subset\II$. If one of
  $p$ or $q$ is an\/ $\II$-primitive icosian, the other must
  be in\/ $\II$ as well.
\end{fact}
\begin{proof}
This follows from \cite[Prop.~1 and Remark~1]{BM2}, where this is
shown for any maximal order of class number one. In particular, it
applies to $\II$.
\end{proof}

As a result of independent interest, we note the following.
\begin{fact} \label{char-poly}
   The linear map $T$ defined by $x\mapsto \alpha\ts px\widetilde{p}$
   has trace $\alpha\ts \N (\tr(p))$ and determinant $\alpha^4\ts\N
   \bigl(\lvert p \rvert^4 \bigr)$, and its characteristic polynomial
   reads
\[ 
   X^4 - \mathrm{trace} (T) X^3 + A X^2 - B X + \det (T)
\]
   where $A=\alpha^2 \bigl( \mathrm{Tr} \bigl( (\tr (p))^2 (\nr (p))'
   \bigr) - 2\ts \N(\nr (p))\bigr)$ and $B= \alpha^3\ts \N \bigl(\tr (p)
   \nr (p)\bigr)$.
\end{fact}
\begin{proof}
This is a straight-forward calculation, e.g., with the matrix
representation from \cite{KR,BM2}, taking into account that
$\nr (\widetilde{p}) = (\nr (p))'$ and expressing the coefficients
in terms of traces and norms.
\end{proof}

In view of Corollary~\ref{para-2}, we now need to understand how an
SSL of $\alat$ of the form $p \alat \widetilde{p}$ with an
$\II$-primitive quaternion relates to the primitive sublattices of
$\alat$. Recall that $\gL$ is a sublattice of $\gG$ if and only if
there is a non-singular integer matrix $Z$ that maps a basis of $\gG$
to a basis of $\gL$, written as $B_{\gL} = B_{\gG} Z$ in terms of
basis matrices, compare \cite{Cassels}. The index is then
\[
    [ \gG : \gL ] \; = \; \lvert \det(Z)\rvert
\]
which does not depend on the actual choice of the lattice basis.  A
sublattice $\gL$ of $\gG\subset\RR^d$ is called
\emph{$\gG\!$-primitive} (or simply primitive) when
\begin{equation} \label{lattice-content}
   \mathrm{cont}^{}_{\gG} (\gL) \; := \;
   \lcm\ts \{ m\in\NN\mid \gL\subset m\gG \} \; = \; 1\, .
\end{equation}
This is equivalent to saying that $\gcd\ts (Z) :=
\gcd\ts \{ Z_{ij}\mid 1\le i,j\le d \} = 1$.

\begin{prop}  \label{prim-one}
  If $p\in\II$ is\/ $\II$-primitive, $p\alat\widetilde{p}$ is an 
  $\alat$-primitive sublattice of $\alat$.
\end{prop}
\begin{proof}
By Lemma~\ref{all-do}, we know that $p\alat\widetilde{p}
\subset\alat$. We thus have to show that $\frac{1}{m}
p\alat\widetilde{p} \subset\alat$ implies $m=1$.  Note first that
$\frac{1}{m}\, p\alat\widetilde{p} \subset\alat$ also implies
$\frac{1}{m}\, p\alat[\tau]\widetilde{p} \subset\alat[\tau]$, with
$\alat[\tau] = \alat \op \tau\alat$ as before. From
\eqref{taulat-in-ico}, we know that $5\ts\II\subset\alat[\tau]$, so
that we get
\[
   \tfrac{5}{m}\, p\ts\II\widetilde{p} \, \subset \,
   \tfrac{1}{m}\, p\alat[\tau]\widetilde{p} \, \subset \,
   \alat[\tau] \, \subset \, \II \, .
\]
Since $p$ is $\II$-primitive by assumption, then so is
$\widetilde{p}$.  By Fact~\ref{no-escape}, this forces $5/m$ to be an
element of $\ZZ[\tau]$.  With $m\in\NN$, this only leaves $m=1$ or $m=5$.

Observe next that $2\alat[\tau]\subset\LL=\ZZ[\tau]^4$. On the other hand,
it is easy to check explicitly (e.g., by means of the basis vectors)
that $\sqrt{5}\,\LL \subset \alat[\tau]$. Together with $2\ts\II\subset
\LL$ from Fact~\ref{L-in-I}, this gives
\[
   \tfrac{4\sqrt{5}}{m}\, p\ts\II\widetilde{p} \, \subset \,
   \tfrac{2\sqrt{5}}{m}\, p\LL\widetilde{p} \, \subset \,
   \tfrac{2}{m}\, p\alat[\tau]\widetilde{p} \, \subset \,
   2\ts \alat[\tau] \, \subset \, \LL \, \subset \, \II \, .
\]
Invoking Fact~\ref{no-escape} again, we see that $\tfrac{4\sqrt{5}}{m}
\in\ZZ[\tau]$, which (with $m\in\NN$) is only possible for $m|4$. In
combination with the previous restriction, this implies $m=1$.
\end{proof}

Combining the results of Corollary~\ref{para-2} and
Proposition~\ref{prim-one}, we obtain the following important
observation.
\begin{coro} \label{prim-two}
   The $\alat$-primitive SSLs of $\alat$ are precisely the ones of the
   form $p\alat\widetilde{p}$ with $p$ an $\II$-primitive
   icosian.  
\end{coro}
\begin{proof}
After Proposition~\ref{prim-one}, it remains to show that every primitive 
SSL $M$ of $\alat$ is of the form $p L \widetilde{p}$ for some primitive 
$p \in \II$. By Corollary~\ref{para-2},  $M = \alpha p L \widetilde{p}$ 
with $\alpha\in\QQ$ and $p\in\II$ primitive. By Proposition~\ref{prim-one}, 
$p  L \widetilde{p}$ is already a primitive sublattice, so $\alpha = \pm 1$.
\end{proof}

Next, we need to find a suitable bijection that permits us to count the
primitive SSLs of $\alat$ of a given index.  Recall from \cite{MW,CMP}
that the unit group of $\II$ has the form
\[
   \II^{\uu} \, = \, \{ x\in\II \mid \N (\nr (x)) = \pm 1 \}
   \, = \, \{ \pm \tau^m \ts \varepsilon \mid m\in\ZZ \mbox{ and }
   \varepsilon\in \varDelta_{H_4} \} \, ,
\]
where $\varDelta_{H_4}$ is defined as above via Eq.~\eqref{h4-vectors}
and satisfies $\varDelta_{H_4} = \II^{\uu}\cap \mathbb{S}^3 = \II \cap
\mathbb{S}^3$, see \cite{BM2,MW} for details\footnote{Let us take this
  opportunity to mention that \cite{BM2}, in deviation from other
  conventions and from the one adopted here, uses the notation
  $\II^{\uu}$ for the $\II$-units on $\mathbb{S}^3$, i.e., for the
  elements of $\varDelta_{H_4}$ only.}.

\begin{lemma} \label{only-units-do}
  For $p\in\II$, one has $p\alat\widetilde{p} = \alat$
  if and only if $p\in\II^{\uu}$.
\end{lemma}
\begin{proof}
  We have $p\alat\widetilde{p}\subset\alat$ for all $p\in\II$ by
  Lemma~\ref{all-do}. The corresponding index is given by the
  determinant of the mapping $x\mapsto px\widetilde{p}$. By
  Fact~\ref{char-poly}, the equality $\alat = p\alat\widetilde{p}$ is
  then equivalent to $\N\bigl(\lvert p\rvert^4\bigr)=1$ and hence to
  $\N (\nr (p)) = \pm 1$, which, in turn, is equivalent to
  $p\in\II^{\uu}$.
\end{proof}

We need one further result to construct a bijective correspondence
between primitive SSLs of $\alat$ and certain right ideals of $\II$, which 
will then solve our problem.

\begin{lemma} \label{symmetry-fits}
   For primitive $r,s \in \II$, one has $r\II = s\II$ if and only if
   $r\alat\ts\widetilde{r} = s\alat\ts\widetilde{s}$.
\end{lemma}
\begin{proof}
  Since $r,s\in\II$, it is clear that $r\II = s\II \;\Rightarrow\; \II
  = r^{-1}s\II\; \Rightarrow\; r^{-1}s \in\II$. Similarly, $s^{-1}r
  \in \II$, so $ r^{-1}s \in \II^\times$.  Lemma~\ref{only-units-do}
  now implies $r^{-1}s \ts L \ts (r^{-1}s)^{\up} = L$, which gives $ r
  L\widetilde{r} = s L \widetilde{s} $.

  In the reverse direction, suppose that $r L\widetilde{r} = s L
  \widetilde{s}$, which gives $y L \widetilde{y} = L$ with $y:=
  r^{-1}s$. Choose $\alpha \in K$ so that $\alpha y$ is a primitive
  element of $\II$, which is always possible. Then,
\[
    \alpha y L \ts\widetilde{\alpha y} \, = \, 
    \alpha \alpha^{\ts\prime} y L \widetilde{y} 
    \, = \, \alpha \alpha^{\ts\prime} L 
\]
is a primitive sublattice of $\alat$ by Proposition~\ref{prim-one}, whence 
$\alpha \alpha^{\ts\prime} = \pm 1$ and $\alpha y L
\ts \widetilde{\alpha y} = L$.

This implies $\varepsilon:= \alpha y =\alpha r^{-1} s \in \II^{\times}$
by Lemma~\ref{only-units-do}. Then,
\[
    r \, = \, r \varepsilon \varepsilon^{-1} \, = \, 
    \alpha s \varepsilon^{-1}  
   \, \in \, \alpha \ts\II \ts ,
\]
so that $\alpha^{\ts\prime} r \in \II$ due to
$\alpha\alpha^{\ts\prime}=\pm 1$, where $\alpha^{\ts\prime}\in K$ by
construction. Since $r\in\II$ is primitive as well, such a relation is
only possible with $\alpha \in \ZZ[\tau]$, in view of the properties
of the $\II$-content of $r$.  Consequently, $\alpha \alpha^{\ts\prime}
= \pm 1$ now gives $\alpha\in\ZZ[\tau]^{\times}$, so that $y$ is an
element of $\II$. Lemma~\ref{only-units-do} now implies $y\in
\II^\times$, whence $y\II = \II$, and finally $s\II = r\II$.
\end{proof}   

In view of our discussion so far, we call a right ideal of $\II$
\emph{primitive} if it is of the form $p\II$ for some primitive
$p\in\II$.
\begin{prop}\label{fundamental-prop}
    There is a bijective correspondence between the primitive right
    ideals of $\II$ and the primitive SSLs of $\alat$, defined by
    $ p \II \mapsto p \alat \widetilde{p} $. Furthermore,
    one has the index formula
\[
   [\II : p\ts\II]  \, = \, \N \bigl( \nr (p)^2 \bigr)
   \, = \, \N \bigl( \lvert p \rvert^4 \bigr) 
    \, = \, [\alat : p\alat\widetilde{p}\,]
\]   
    under this correspondence.
\end{prop}
\begin{proof}
  It is clear from Lemma~\ref{symmetry-fits} that the mapping is
  well-defined and injective, while Corollary~\ref{prim-two} implies
  its surjectivity.

The index of $p\ts\II$ in $\II$ derives from the determinant
formula in Fact~\ref{quat-simi} by taking the norm in
$\ZZ[\tau]$.  The index relation now follows from Fact~\ref{char-poly}.
\end{proof}

\begin{remark} \label{symmetries} 
  The results of Lemmas~\ref{only-units-do} and \ref{symmetry-fits}
  also show that $\varDelta_{H_4} = \II^{\times}\cap \, \mathbb{S}^3$,
  which is a subgroup of $\II^{\times}$ of order $120$, can be viewed
  as the standard double cover of the rotation symmetry group of
  $\alat$, which has order $60$ and contains the point reflection in
  the origin (note that $-1$ is a rotation in $4$-space). This is a
  geometric interpretation of the bijective correspondence.  Let us,
  in this context, also point out that $p\II\widetilde{p}\cap V_{+} =
  p\alat\widetilde{p}$, but that $p\II\cap V_{+}$ is generally a much
  bigger set.
\end{remark}

With the representation of a general SSL as an integer multiple of a
primitive SSL mentioned earlier, and observing that all possible
indices of SSLs are squares, we can now solve our original problem.
We need the Dirichlet character
\begin{equation} \label{character}
   \chi (n) \; = \;
   \begin{cases} 0, &  n \equiv 0 \; (5)\, , \\
                 1, &  n \equiv \pm 1 \; (5) \, , \\
                -1, &  n \equiv \pm 2 \; (5) \, . 
   \end{cases}
\end{equation}
Note that the corresponding $L$-series, $L(s,\chi) = \sum_{m=1}^{\infty}
\,\chi(m)\ts m^{-s}$, defines an entire function on the complex plane.
Furthermore, we need the zeta functions of the icosian ring,
\begin{equation} \label{ico-zetas}
    \zeta^{}_{\II} (s) \; = \; \zeta^{}_{K} (2s) \, \zeta^{}_{K} (2s-1)
    \quad \mbox{and} \quad \zeta^{}_{\II.\II} (s) \; = \;
    \zeta^{}_{K} (4s) \, ,
\end{equation}
where $\zeta^{}_{\II} (s)$ and $\zeta^{}_{\II.\II} (s)$ denote the
Dirchlet series for the one-sided and the two-sided ideals of $\II$,
respectively. Here, $\zeta^{}_{K} (s) = \zeta(s)\, L(s,\chi)$ is the
Dedekind zeta function of the quadratic field $K$, see
\cite{V,BM2,BPR} for details.

\begin{theorem}
    The number of SSLs of a given index is the same for the lattices
    $A_4$ and $\alat$. The possible indices are the squares of
    non-zero integers of the form $k^2 + k\ell - \ell^2 =
    \N(k+\ell\tau)$. Moreover, all possible indices are realized.

    {}Furthermore, there is an index preserving bijection between the
    primitive SSLs of $A_4$ and the primitive right ideals of\/ $\II$.
    When $f(m)$ denotes the number of SSLs of index $m^2$ and
    $f^{\sf pr} (m)$ the number of primitive ones,
    the corresponding Dirichlet series generating functions read
\[
      D_{A_4} (s) \; := \; \sum_{m=1}^{\infty} \frac{f(m)}{m^{2s}}
      \; = \; \zeta(4s)\, \frac{\zeta^{}_{\II}(s)}{\zeta^{}_{K} (4s)}
      \; = \; \frac{\zeta^{}_{K} (2s)\,\zeta^{}_{K} (2s-1)}{L (4s,\chi)} 
\]  
and
\[
      D_{A_4}^{\sf pr} (s) \; := \; 
     \sum_{m=1}^{\infty} \frac{f^{\sf pr}(m)}{m^{2s}}
      \; = \; \zeta^{\sf pr}_{\II} (s)
     \; = \; \frac{\zeta^{}_{\II} (s)}{\zeta^{}_{K} (4s)}
     \; = \; \frac{\zeta^{}_{K} (2s) \, 
     \zeta^{}_{K} (2s-1)}{\zeta^{}_{K} (4s)} ,
\]
   with $\zeta^{}_{K} (s)$ and $L(s,\chi)$ as defined above.
   In particular, each possible index is thus also realized
   by a primitive SSL.
\end{theorem}   

\begin{proof}
  The first claim follows from Fact~\ref{invar-similar}, while the
  index characterization follows either from \cite{CRS} or from the
  index formula in Proposition~\ref{fundamental-prop}.  Recall that
  $k^2 + k\ell - \ell^2 = \N (k+\ell\tau)$ is the norm the principal
  ideal $(k+\ell\tau)\ts\ZZ[\tau]$, and that all $\ZZ[\tau]$-ideals
  are of this form.  Since the reduced norm as a mapping from icosian
  right ideals to ideals in $\ZZ[\tau]$ is surjective, each possible
  index is realized.
  
  In view of the previous remarks, we have
\[
    D^{}_{A_4} (s) \; = \; \zeta(4s)\, D^{\sf pr}_{A_4} (s) \, ,
\]
where, due to our bijective correspondence and the index formula of
Proposition~\ref{fundamental-prop}, the various formulas for the
Dirichlet series are immediate, see \cite{BM2,BPR} for a derivation of
$\zeta^{\sf pr}_{\II} (s)$ as the quotient stated above. The last
equality then follows from \eqref{ico-zetas}. Since $f^{\sf pr}(m)$
vanishes precisely when $f (m)$ does (see below for an explicit
formula), the last claim is clear.
\end{proof}

Inserting the Euler products of $\zeta(s)$ and $\zeta^{}_{K} (s)$,
one finds the expansion of the Dirichlet series $D_{A_4} (s)$ 
as an Euler product,
\[
   D_{A_4} (s) \; = \; \frac{1}{(1-5^{-2s})(1-5^{1-2s})}
   \prod_{p\equiv\pm 1\; (5)} \frac{1+p^{-2s}}{1-p^{-2s}}\,
   \frac{1}{(1-p^{1-2s})^2} \prod_{p\equiv\pm 2\; (5)}
   \frac{1+p^{-4s}}{1-p^{-4s}}\,  \frac{1}{1-p^{2-4s}} 
\]
and similarly for $D_{A_4}^{\sf pr} (s)$. Consequently, the arithmetic 
function $f(m)$ (and also $f^{\sf pr} (m)$) is multiplicative, i.e.,
$f(mn) = f(m)\, f(n)$ for $m,n$ coprime, with $f(1)=1$. The function
$f$ is then completely specified by its values at prime powers $p^r$
with $r\ge 1$. These are given by 
\[
    f(p^r) \; = \; \begin{cases}
        \frac{5^{r+1} - 1}{4} , & \text{if $p=5$,} \\[1mm]
      \frac{2\, (1-p^{r+1}) - (r+1)(1-p^2)p^r}{(1-p)^2}, &
      \text{for primes $p\equiv\pm 1 \; (5)$,} \\[1mm]
      \frac{2 - p^{r} - p^{r+2}}{1-p^2} , &
      \text{for primes $p\equiv\pm 2 \; (5)$ and $r$ even,} \\[1mm]
      0, &
      \text{for primes $p\equiv\pm 2 \; (5)$ and $r$ odd,}
    \end{cases}
\]
and similarly for $f^{\sf pr} (m)$:
\[
   f^{\sf pr} (p^r) \; = \; \begin{cases}
        6\cdot 5^{r-1}, & \text{if $p=5$,} \\
       (r+1)p^r+2rp^{r-1}+(r-1)p^{r-2}, & \text{for $p \equiv \pm 1\; (5)$,} \\
       p^r+p^{r-2}, & \text{for $p \equiv \pm 2 \;(5)$ and $r$ even,} \\
       0, & \text{for $p \equiv \pm 2 \;(5)$ and $r$ odd.}
\end{cases}
\]
The first few terms of the Dirichlet series thus read
\[
\begin{split}
   D_{A_4} (s) & \; = \; 1 + \tfrac{6}{4^{2s}} + \tfrac{6}{5^{2s}} + 
   \tfrac{11}{9^{2s}} + \tfrac{24}{11^{2s}} + \tfrac{26}{16^{2s}} + 
   \tfrac{40}{19^{2s}} + \tfrac{36}{20^{2s}} + \tfrac{31}{25^{2s}} +
   \tfrac{60}{29^{2s}} + \tfrac{64}{31^{2s}} + \tfrac{66}{36^{2s}}
   + \ldots \\
   D_{A_4}^{\sf pr} (s) & \; = \; 1 + \tfrac{5}{4^{2s}} + \tfrac{6}{5^{2s}} + 
   \tfrac{10}{9^{2s}} + \tfrac{24}{11^{2s}} + \tfrac{20}{16^{2s}} + 
   \tfrac{40}{19^{2s}} + \tfrac{30}{20^{2s}} + \tfrac{30}{25^{2s}} +
   \tfrac{60}{29^{2s}} + \tfrac{64}{31^{2s}} + \tfrac{50}{36^{2s}}
   + \ldots
\end{split}   
\]
where the denominators are the squares of the integers previously
identified in \cite{CRS}, see also entry {\sf A{\ts}031363} of
\cite{online}.  Further details can be found in \cite{H}.

These series may now be compared with the zeta functions
$\zeta^{}_{\II}(s)$ and $\zeta^{\sf pr}_{\II}(s)$, which reveals the
origin of the various contributions. In particular, the $6$ SSLs of
index $16$ stem from the $5$ generators of primitive ideals of $\II$
of index $16$ together with the SSL $2 A_4$. No such extra solution
exists for index $25$, while index $81$ emerges also from the SSL $3 A_4$.

Due to our definition with $f(m)$ being the number of SSLs of $A_4$ of
index $m^2$, the Dirichlet series generating function of the
arithmetic function $f$ is $D_{A_4} (s/2)$, which has nice analytic
properties. In particular, it is analytic on the half-plane
$\{\sigma>2\}$, where we write $s = \sigma + \ii t$ as usual.
Moreover, it is analytic on the line $\{\sigma = 2\}$, except at
$s=2$, where we have a simple pole as the right-most singularity of
$D_{A_4} (s/2)$.  Consequently, one can derive the asymptotic growth
of $f(m)$ from it.  Since the value of the arithmetic function $f(m)$
fluctuates heavily, this is done via the corresponding summatory
function
\[
      F(x) \; := \; \sum_{m\le x}  f(m) \; \sim \; 
      \varrho\, \frac{x^2}{2}\, ,   \quad\text{as $x\to\infty$},
\]
where the growth constant is given by
\[
   \varrho \; = \; \mathrm{res}_{s=2} \, D_{A_4} (s/2)
   \; = \; \frac{\zeta^{}_{K} (2)\, L(1,\chi)}{L(4,\chi)}
   \; = \; \frac{1}{2}\sqrt{5}\,\log(\tau)
   \; \simeq \; 0.538\, 011 .
\]
This result relies on Delange's theorem, see \cite[Appendix]{BM2} 
for details.  The corresponding calculations for the asymptotic behaviour
of $f^{\sf pr} (m)$ are analogous.

\section{Results for related lattices}

{}From previously published results, one can read off or easily
derive the generating functions for the root lattices
$A_d$ with $d\le 3$.
\begin{theorem}
   The Dirichlet series generating functions for the 
   number of SSLs of the root lattices
   $A_d$ with $d\le 3$ are $D^{}_{\!A_1} (s) = \zeta(s)$,
   $D^{}_{\!A_2} (s) = \zeta^{}_{\QQ(\xi_3)}$, with $\xi_3=e^{2\pi i/3}$
   and $\zeta^{}_{\QQ(\xi_3)} (s)$ the Dedekind zeta function
   of the cyclotomic field $\QQ(\xi_3)$, and
   $ D^{}_{\!A_3} (s) = \zeta(3s) \, \varPhi_{\rm cub} (3s)$,
   where 
\[
    \varPhi_{\rm cub} (s) \; = \; 
    \frac{1-2^{1-s}}{1+2^{-s}}\;
    \frac{\zeta(s)\,\zeta(s-1)}{\zeta(2s)} 
\]
   is the generating function of the related cubic coincidence 
   site lattice problem derived in \cite{B,BPR}.
\end{theorem}

\begin{proof}
  $A_1$ is a scaled version of the integer lattice $\ZZ$, whence
  $D_{\!A_1} (s) = D_{\ZZ} (s) = \zeta(s)$ follows from
  Fact~\ref{invar-similar} together with \cite[Prop.~3.1]{BM1}. The
  triangular lattice $A_2$ is a scaled version of the ring of
  Eisenstein integers in $\QQ(\xi_3)$, so that this generating
  function follows from \cite[Prop.~1]{BG}.  It can also be written as
  a product of two Dirichlet series,
\[
    D_{\!A_2} (s) \; = \; \zeta(2s)\,
    \frac{\zeta^{}_{\QQ(\xi_3)} (s)}{\zeta(2s)} ,
\]
one for the scaling by integers and the other for the primitive SSLs.

{}For the primitive cubic lattice $\ZZ^3$, the SSL generating function
was previously identified as $\zeta(3s)\, \varPhi_{\rm cub}(3s)$ in
\cite[Thm.~5.1]{BM1}. Noting that $A_3$ is a version of the
face-centred cubic lattice, compare \cite{CS}, one has to check that
the cubic lattices share the same SSL statistics. This is a
straight-forward calculation with their basis matrices and the
integrality conditions for similarity transforms, similar to the one
outlined in \cite[Sec.~5]{BM1}, using the Cayley parametrization of
matrices in $\mathrm{SO} (3,\QQ)$ from \cite{B}, see also \cite{H} for
details.
\end{proof}

Clearly, by Fact~\ref{invar-dual}, one then also has the relation
\begin{equation}
   D_{\!A^{}_3} (s) \; = \; D_{\!A^*_3} (s) \; = \;
   D_{\ZZ^3} (s) \, .
\end{equation}
This identity can also be derived \cite{Z} by a consideration of the
shelling structure of $A_{3}^{*}$ relative to $\ZZ^3$, which bypasses
the explicit basis calculations referred to above.

Other completely worked out examples include the square lattice
$\ZZ^2$, with $D_{\ZZ^2} (s) = \zeta^{}_{\QQ(i)} (s)$, see \cite{BM1},
and the cubic lattices in $d=4$, see \cite{BM2}.  Also, various
$\ZZ$-modules of rank $r>d$ in dimensions $d\le 4$ are solved in
\cite{BM1,BM2,BM3,BG}. Common to all these examples is the rather
explicit use of methods from algebraic number theory or quaternions in
conjunction with suitable parametrizations of the rotations involved.
In higher dimensions, results are sparse, compare \cite{CRS} and
references given there, and we are not aware of any complete solution
in terms of generating functions at present.

\section{Algebraic afterthoughts} \label{after}

The special twist map $\ts\widetilde{.}\ts$ played a crucial role in
solving the sublattice problem above. It is an odd looking map at
first sight, but has a number of interesting properties that are worth
considering in more detail. To do so, we also write $\eta (x) =
\widetilde{x}$ from now on, and use both notations in parallel.
Moreover, we define inner automorphisms $T_a$ of the quaternion
algebra $\HH (K)$ via $T_a (x) = axa^{-1}$, for $a\in\HH (K)^\bt$.
Finally, we call any $K$-semilinear involutory anti-automorphism of
$\HH (K)$ a \emph{twist map}.
\begin{lemma} \label{twist-products}
  If $\sigma^{}_{1}$ and $\sigma^{}_{2}$ are two twist maps of\/
  $\HH(K)$ that map $\II$ into itself, their product is an inner
  automorphism of\/ $\II$, i.e., $\sigma^{}_{1}\ts \sigma^{}_{2}
  = T_{\varepsilon}$ for some $\varepsilon\in\II^{\uu}$. In particular,
  any such twist map is of the form $T_{\varepsilon} \eta$ for some
  $\varepsilon\in\II^{\uu}$.
\end{lemma}
\begin{proof}
  Since both $\sigma_i$ are $K$-semilinear, their product is
  $K$-linear. Also, the product $\sigma^{}_{1} \sigma^{}_{2}$ of two
  anti-automorphisms is an automorphism of $\HH (K)$, which must be
  inner by the Skolem-Noether theorem \cite[Thm.~I.1.4]{Rost}. This
  gives $\sigma^{}_{1}\ts \sigma^{}_{2} = T_a$ with $a\in\HH (K)^\bt$.
  By construction, $T_a$ is also an automorphism of $\II$, whence we
  may assume $a\in\II$ without loss of generality. This implies $a\II
  = \II a$, so that $a\II$ is a two-sided ideal. Within $\II$, this
  implies $a=\alpha\ts\varepsilon$ with $0\neq\alpha\in\ZZ[\tau]$ and
  $\varepsilon\in\II^{\uu}$. Observing $T_{\alpha\varepsilon}=
  T_{\varepsilon}$ establishes the claim.
\end{proof}

For our further arguments, it is advantageous to consider the mapping
$\vartheta \! : \; \HH(K)\longrightarrow V_{+}$, defined by $x\mapsto
x\,\eta(x)=x\widetilde{x}$, where $V_{+}=\{x\in\HH(K)\mid\widetilde{x}=x\}$
as before. Clearly, one has $\vartheta(\II)\subset\alat$, as a consequence 
of Proposition~\ref{key-prop}. 
\begin{lemma} \label{restricted-map} 
  The mapping $\vartheta$ satisfies
  $\vartheta(\II^{\uu}) = \vartheta(\varDelta_{H_4}) \subset
  \varDelta_{A_4}$. Moreover, the corresponding restriction,
  $\vartheta_{\varDelta}\!:\,
  \varDelta_{H_4}\longrightarrow\varDelta_{A_4}$, is a surjective
  $6$-to-$1$ mapping.
\end{lemma}
\begin{proof}
  Any $a\in\II^{\uu}$ has the form $a=\pm \tau^m \varepsilon$ for some
  $m\in\ZZ$ and $\varepsilon\in\varDelta_{H_4}$, so that
  $a\widetilde{a} = (\tau \tau')^m \, \varepsilon\tilde{\varepsilon} =
  (-1)^m \, \varepsilon\tilde{\varepsilon}$. The identity
  $\vartheta(\II^{\uu}) = \vartheta(\varDelta_{H_4})$ now follows from
  the observation that $\ii \in \varDelta_{H_4}$, with $\vartheta(\ii)
  = \ii^2 = -1$. If $\varepsilon\in\varDelta_{H_4}$,
  $\varepsilon\tilde{\varepsilon}$ lies both in $\alat$ and in
  $\varDelta_{H_4}$, and thus in $\varDelta_{A_4}=\alat \cap
  \varDelta_{H_4}$, see Eq.~\eqref{root-systems}.

  Observe next that the equation $\varepsilon\tilde{\varepsilon} = 1$
  has precisely $6$ solutions with $\varepsilon\in\varDelta_{H_4}$,
  namely $\pm 1$ together with $\pm\tfrac{1}{2} (\pm
  1,0,\tau,\tau\!-\!1)$. Under multiplication, they form the cyclic
  group $C_6$, generated by the quaternion
  $\tfrac{1}{2}(1,0,\tau,\tau\!-\!1)$.  Clearly, one has
  $\vartheta(a)=\vartheta(b)$ with $a,b\in\varDelta_{H_4}$ if and only
  if $b=a\varepsilon$ with $\varepsilon\in C_6$, so that a coset
  decomposition of the $120$ roots of $\varDelta_{H_4}$ with respect
  to $C_6$ establishes the second claim.
\end{proof}   
 
Define now, for an arbitrary $\varepsilon\in\II^{\uu}$, the mapping
$\eta_{\varepsilon} = T_{\varepsilon}\eta T_{\varepsilon^{-1}} $.
Since $\eta_{\gamma\varepsilon} = \eta_{\varepsilon}$ for any
$\gamma\in\ZZ[\tau]^{\uu}$, all such maps can be obtained from an
$\varepsilon\in\varDelta_{H_4}$.

\begin{prop} \label{all-twists} 
  The twist maps of\/ $\HH (K)$ that
  map\/ $\II$ into itself are precisely the maps
  $\eta_{\varepsilon}=T_{\varepsilon}\eta T_{\varepsilon^{-1}} =
  T^{}_{\varepsilon\tilde{\varepsilon}}\, \eta $, with
  $\varepsilon\in\varDelta_{H_4}$.  Equivalently, they are precisely
  the maps $T_a \eta$ with $a\in\varDelta_{A_4}$.
\end{prop}
\begin{proof}
  By Lemma~\ref{twist-products}, we know that any twist map is of the
  form $T_a \ts\eta$ for some $a\in\II^{\uu}$. It is easy to check the
  commutation relation $\eta\ts T^{}_{a^{-1}} =
  T^{}_{\widetilde{a}}\ts\eta$, which also implies, via left
  multiplication by $T_a$, that
\[
    T^{}_{a}\eta\ts T^{}_{a^{-1}} \; = \; T^{}_{a\widetilde{a}}\ts\eta\, . 
\]
A twist map is an involution (in particular, one has
$\eta^2=\mathrm{Id}$), so that $ (T_a \eta)^2 = T_{a
  (\widetilde{a})^{-1}}$ has to be the identity.  Since $T_b =
\mathrm{Id} = T_1$ with $b\in\II^{\uu}$ implies $b\in\II^{\uu}\cap K
=\ZZ[\tau]^{\uu}$, the previous condition means $\widetilde{a} = \beta
a$ for some $\beta\in\ZZ[\tau]^{\uu}$. Consequently, $a= (\beta
a)^{\up} = \beta'\, \widetilde{a} = \beta' \beta a$, whence $\beta'
\beta = 1$. This implies $\beta = \pm\tau^{2m}$ for some $m\in\ZZ$.
Observing that $T_{\tau^m a} = T_{a}$, we may thus, without loss of
generality, replace $a$ by $\tau^m a$, so that the involution
condition reduces to $\widetilde{a}=\pm a$.

Next, write $a =\gamma \varepsilon$ with $\gamma \in\ZZ[\tau]^{\uu}$
and $\varepsilon\in\varDelta_{H_4}$, which is possible due to the
structure of $\II^{\uu}$. Inserting this into $\widetilde{a}=\pm a$
results in $\gamma^{-1} \gamma' = \pm \varepsilon
\tilde{\varepsilon}^{-1}$. Observe next that $\ZZ[\tau]^{\uu}$ and
$\varDelta_{H_4}$ share only $\pm 1$. Within $\ZZ[\tau]^{\uu}$,
$\gamma^{-1} \gamma' = - 1$ has no solution, while
$\widetilde{\varepsilon}=-\varepsilon$ has no solution in
$\varDelta_{H_4}$, meaning that we must have
$\gamma^{-1} \gamma' = 1 = \varepsilon
\tilde{\varepsilon}^{-1}$.
This leaves us with the condition $\widetilde{a}=a$, which is
precisely solved by the $20$ elements of $\varDelta_{A_4}$, see
Eq.~\eqref{root-systems}. All possible twist maps are thus of the form
$T_a \eta$ with $a\in\varDelta_{A_4}$.

Observing $\varepsilon\tilde{\varepsilon}=\vartheta(\varepsilon)$,
the first claim now follows from Lemma~\ref{restricted-map}.
\end{proof}

\begin{theorem} \label{ten-lat} 
  There are exactly $10$ distinct twist
  maps of\/ $\HH(K)$ that map $\II$ into itself, which may be written
  as $\eta_{\varepsilon}$ with $\varepsilon$ running through the $10$
  positive roots of $\varDelta_{A_4}$.
  
  Moreover, the fixed points of $\eta_{\varepsilon}$ are precisely the
  points of the lattice $T_{\varepsilon} (\alat)$, each of which is a
  root lattice of type $A_4$.
\end{theorem}
\begin{proof}
  The twist maps $T_a \eta$ with $a\in\varDelta_{A_4}$ from
  Proposition~\ref{all-twists} form precisely $10$ distinct maps
  because $T_{b} = T_{a}$ holds here if and only if $b=\pm a$.
  Consequently, the restriction to the positive roots of
  $\varDelta_{A_4}$ suffices.
  
  From Proposition~\ref{key-prop}, we know that $\alat =
  \{x\in\II\mid\widetilde{x}=x\}$ is the lattice of fixed points of
  $\eta^{}_{1} = \eta$. It is an immediate consequence that
  $\eta_{\varepsilon}$ fixes the points of
  $\varepsilon\alat\varepsilon^{-1} = T_{\varepsilon} (\alat)$.
\end{proof}

\begin{remark}
  There are $60$ congruent copies of $L$ within the icosian ring
  $\II$, each of which must be of the form $aLb$ with 
  $a,b\in\varDelta_{H_4}$ due to the structure of $\II^{\times}$. 
  Observing $a\alat\widetilde{a}=\alat= -\alat=-a\alat\widetilde{a}$
  for $a\in\varDelta_{H_4}$, it is evident that the
  stabilizer of $\alat$, when seen as a subgroup of 
  $\varDelta_{H_4}\!\times\!\varDelta_{H_4}$, has order $240$.
  Since the order of $\varDelta_{H_4}\!\times\!\varDelta_{H_4}$
  is $120^2$, the length of the orbit of $\alat$ under the action
  of this group is $60$, which are the congruent copies.
  Among them, the $10$ lattices of Theorem~\ref{ten-lat} are
  precisely the ones that contain $1$. They are thus special in the
  sense that they line up with the arithmetic of $\II$.
\end{remark}

Consider $T\! : \varDelta_{H_4} \longrightarrow \mathrm{Aut} (\II)$,
$\varepsilon\mapsto T_\varepsilon$, which is a group homomorphism. As
$T_{\varepsilon} = T_{-\varepsilon}$, $T(\varDelta_{H_4})\simeq Y$ is
the standard icosahedral group of order $60$. If $c$ denotes the
conjugation map, $x\mapsto c(x):=\bar{x}$, one has $c\,
T_{\varepsilon} = T_{\varepsilon}\ts c$ for all
$\varepsilon\in\varDelta_{H_4}$, so that
\[
    \bigl\langle T(\varDelta_{H_4}), c \,\bigr\rangle \; \simeq \;
    Y\times C_2 \; = \; Y_h
\]
is the full icosahedral group of order $120$. By Lemma~\ref{props}, $c$
and $\eta$ commute as well, while $\eta T_{\varepsilon}\eta =
T_{\tilde{\bar{\varepsilon}}}$ shows that
\[
    \bigl\langle T(\varDelta_{H_4}), c, \eta\, \bigr\rangle \; \simeq \;
    Y_h \rtimes \langle \eta\rangle \; = \;
    \bigl( Y\rtimes \langle\eta\rangle \bigr)\times\langle c\rangle\, .
\]

\begin{lemma} 
 One has $Y\rtimes \langle\eta\rangle \simeq S_5$, with
  $S_5$ being the full permutation group of\/ $5$ elements.
\end{lemma}
\begin{proof}
  Note that $Y$ is isomorphic with the alternating group of $5$
  elements, which is simple.  As $\langle\eta\rangle\simeq C_2$,
  $Y\rtimes \langle\eta\rangle$ is a $C_2$-extension of $Y$, where the
  automorphism on $Y$ induced by $\eta$ is $T_{\varepsilon} \mapsto
  T_{\tilde{\bar{\varepsilon}}}$, hence not the identity.  As the
  automorphism class group of $Y$ is the cyclic group $C_2$, there are
  (up to isomorphism) only two $C_2$-extensions of $Y$, namely $Y_h$
  and $S_5$. Since $Y_h\simeq Y\times C_2$, the claim follows.
\end{proof}

Let us add some comments on the geometric meaning of the $10$ twist
maps, in the setting of root systems, where we focus on the interplay
between twist maps that map $\II$ into itself and elements of order
$3$. Here, we use the natural quadratic form $\tr(x\bar{y})$ mentioned
in Remark~\ref{rem-one}.

The intersection of the fixed points of $\eta$ and the pure icosians
$\II_{\ts 0} = \{ x\in\II\mid x+\bar{x}=0 \}$ is the $\ZZ$-span of a root
system of type $A_2$. The latter consists of the $6$ roots
\[
  \pm(0,1,0,0) \, , \; \pm\tfrac{1}{2}(0,\pm1, \tau \!-\! 1, -\tau) \, ,
\] 
which is also the intersection of $\II_{\ts 0}$ with our original
$A_4$ root system from Eq.~\eqref{a4-roots}.

Now, the twist maps are semilinear involutions of $\II$, each of them
a conjugate of our given one, $\eta$, by an element of the group
$\mathbb{I}^{\uu}$.  Since all conjugations by elements of
$\mathbb{I}^\times$ are orthogonal transformations that stabilize $1$
and $\II_{\ts 0}$, each twist map likewise fixes an $A_2$ root system
which lies both inside its own $A_4$ and inside $\II_{\ts 0}$.

We can easily see the geometric meaning of these $10$ root systems of
type $A_2$.  We have noted that each $a\in \varDelta_{H_4}$
determines, by conjugation, an orthogonal transformation $T_a$ that
preserves $\II$, $\II_{\ts 0}$, and $\varDelta_{H_4}$, where
$\varDelta_{H_4} = \II\cap\mathbb{S}^3$ as mentioned above. In
particular, $T_a$ preserves $\II_{\ts 0}\cap \varDelta_{H_4}$, which
is a root system of type $H_3$, with $30$ elements, that we denote by
$\varDelta_{H_3}$. Its convex hull is the icosidodecahedron,
a semi-regular polytope with the $30$ roots as its vertices, $20$
triangular faces, and 12 pentagonal faces.  It is the dual of Kepler's
famous triacontahedron, see \cite[Plate I, Figs.~10 and 12]{Cox} for
details. The $20$ triangular faces, which come in centrally symmetric
pairs, give rise to $10$ axes of $3$-fold symmetry --- which, in fact,
correspond to the elements of order $3$ in the symmetry group $Y_h$ of
the icosahedron.  We shall show that each of these axes has a hexagon
of roots in the plane orthogonal to it. These hexagons are precisely
the root systems of type $A_2$ that arise from our twist maps.

To see this, recall that the elements of the form $T_a$, with $a \in
\varDelta_{H_4 }$, give rise to the symmetry group (of order $120$) of
the root system $\varDelta_{H_3}$, via the mapings $x\mapsto T_a(x)$
and $x\mapsto T_a(\bar{x})$. Here, the first type provides the $60$
orientation preserving elements (isomorphic to the alternating group
on $5$ elements), and the second type the $60$ orientation reversing
elements. Together, they form $Y_h$, the symmetry group of the
icosahedron (and hence also of the icosidodecahedron and the
triacontahedron).

Because all our twist maps are conjugate to one another, it
suffices to look at the $A_2$ lattice indicated above. It has the 
$A_2$ root basis
\[
   \{(0,-1,0,0),\tfrac{1}{2}(0,1, \tau \!-\! 1, -\tau) \}\, ,
\]
and the product of the two reflections generated by these quaternions is 
the rotation of order $3$ that is a Coxeter element of this root system
(the only other one being its inverse, which is obtained by
multiplying the two elements in the opposite order).  It is the inner
automorphism determined by
\[
  -\tfrac{1}{2}(0,1,0,0)(0,1, \tau \!-\! 1, -\tau) \; = \;
  \tfrac{1}{2}(-1, 0, \tau, -\tau') = : z \, . 
\] 
Of course, $z$ is an element of $\varDelta_{H_4}$, and it is easy to see
that $z^3 =-1$, so $T_z$ has order $3$. These inner automorphisms and
their inverses give $20$ elements of order $3$, which then exhaust the
elements of order $3$ in $Y_h$.  This establishes the connection
between the $A_2$ root systems in $\varDelta_{H_3}$ and the rotations of
order $3$ in the icosahedral group.

Each $A_2$ inside $\varDelta_{H_3}$ completes uniquely to an $A_4$
which contains $\pm 1$. This is easy to see from Eq.~\eqref{lat-1}, in
which the four icosians displayed are an $A_4$ basis which includes
the $A_2$ that we have been using in our above argument.
 
So, in short, geometrically we can locate the $A_4$ root systems that
are the fixed points of the special twist maps as follows. Each of the
great circles through $6$ vertices of the icosidodecahedron of roots
in $\II_{\ts 0}$ is a root system of type $A_2$, and it extends
uniquely to a root system of type $A_4$ inside $\varDelta_{H_ 4}$ that
contains $1$.

\smallskip
\section*{Acknowledgements}
It is a pleasure to thank Ulf Rehmann, Rudolf Scharlau, Rainer
Schulze-Pillot and Peter Zeiner for helpful discussions, and Peter
Pleasants for carefully reading the manuscript.  We thank the referee
for a number of helpful suggestions to improve the manuscript. This
work was supported by the German Research Council (DFG), within the
CRC 701, by EPSRC via grant EP/D058465, and by the Natural Sciences
and Engineering Research Council of Canada (NSERC).


\smallskip
\begin{thebibliography}{99}
\small

\bibitem{B}
M.~Baake,
Solution of the coincidence problem in dimensions $d\le 4$,
in:\ \textit{The Mathematics of Long-Range Aperiodic Order},
ed.\ R.\thinspace V.\ Moody, NATO-ASI C 489, Kluwer,
Dordrecht (1997), pp.\ 9--42;
rev.\ version, 
\texttt{arXiv:math.MG/0605222}.

\bibitem{BG}
M.~Baake and U.~Grimm,
Bravais colourings of planar modules with $N$-fold symmetry,
\textit{Z.\ Krist.} {\bf 219} (2004) 72--80;
\texttt{arXiv:math.CO/0301021}.

\bibitem{BM1} 
M.~Baake and R.\thinspace V.~Moody, 
Similarity submodules and semigroups, in:\ \textit{Quasicrystals 
and  Discrete Geometry}, ed.\ J.\ Patera, FIM 10, 
AMS, Providence, RI (1998), pp.\ 1--13.

\bibitem{BM2}
M.~Baake and R.\thinspace V.~Moody,
Similarity submodules and root systems in 
four dimensions,
\textit{Can.\ J.\ Math.} {\bf 51} (1999) 1258--1276;
\texttt{arXiv:math.MG/9904028}.

\bibitem{BM3}
M.~Baake and R.\thinspace V.~Moody,
Invariant submodules and semigroups of self-similarities
for Fibonacci modules, in:\ \textit{Aperiodic {}'$\, 97$}, 
eds.\ M.\ de Boissieu, J.-L.\ Verger-Gaugry
and R.\ Currat, World Scientific, Singapore (1998), pp.\ 21--27;
\texttt{arXiv:math-ph/9809008}.

\bibitem{BPR}
M.~Baake, P.~Pleasants and U.~Rehmann,
Coincidence site modules in $3$-space,
\textit{Discrete Comput.\ Geom.} {\bf 38} (2007) 111--138;
\texttt{arXiv:math.MG/0609793}.

\bibitem{Cassels}
J.\thinspace W.\thinspace S.\ Cassels,
\textit{Introduction to the Geometry of Numbers},
reprint, Springer (1997).

\bibitem{CMP}
L.~Chen, R.\thinspace V.~Moody and J.~Patera,
Non-crystallographic root systems, in:\
\textit{Quasicrystals and Discrete Geometry},
ed.\ J.\ Patera, FIM 10, AMS, Providence, RI (1998),
pp.\ 135--178.

\bibitem{CRS}
J.\thinspace H.~Conway, E.\thinspace M.~Rains and 
N.\thinspace J.\thinspace A.~Sloane,
On the existence of similar sublattices,
\textit{Can.\ J.\ Math.} {\bf 51} (1999) 1300--1306.

\bibitem{CS}
J.\thinspace H.~Conway and N.\thinspace J.\thinspace A.~Sloane,
\textit{Sphere Packings, Lattices and Groups},
3rd ed., Springer, New York (1999).

\bibitem{Cox}
H.\thinspace S.\thinspace M.~Coxeter,
\textit{Regular Polytopes}, 
3rd ed., Dover, New York (1973).

\bibitem{H}
M.~Heuer,
\textit{\"Ahnlichkeitsuntergitter des Wurzelgitters $A_4$},
Diplomarbeit, Univ.\ Bielefeld (2006).

\bibitem{KR} 
M.~Koecher and R.~Remmert, 
Hamilton's quaternions, in:\ \textit{Numbers}, 
eds.\ H.-D.\ Ebbinghaus et al.,
Springer, New York (1991), pp.~189--220.

\bibitem{Rost}
M.-A.~Knus, A.~Merkurjev, M.~Rost and J.-P.~Tignol,
\textit{The Book of Involutions},
Coll.\ Publ.\ 44, AMS, Providence, RI (1998).

\bibitem{MW}
R.~V.~Moody and A.~Weiss,
On shelling $E_8$ quasicrystals,
\textit{J.\ Number Theory}  \textbf{47} (1994) 405--412.

\bibitem{O}
O.~T.~O'Meara,
\textit{Introduction to Quadratic Forms},
corr.\ 3rd printing, Springer, Berlin (1973).

\bibitem{Rehmann}
U.~Rehmann, private communication (2006).

\bibitem{S1}
R.~Scharlau,
A structure theorem for weak buildings of
spherical type,
\textit{Geom.\ Dedicata} \textbf{24} (1987) 77--84.

\bibitem{S2}
R.~Scharlau,
Unimodular lattices over real quadratic fields,
\textit{Math.\ Z.} \textbf{216} (1994) 437--452.

\bibitem{online}
N.\thinspace J.\thinspace A.~Sloane,
\textit{The Online Encyclopedia of Integer Sequences},
\newline
\texttt{http://www.research.att.com/\~{}njas/sequences/}

\bibitem{V}
M.-F.~Vign\'eras, 
\textit{Arithm\'etique des Alg\`ebres de Quaternions},
LNM 800, Springer, Berlin (1980).

\bibitem{Z}
P.~Zeiner, private communication (2006).

\end{thebibliography}
\end{document}